\documentclass[12pt]{amsart}
\usepackage{amsmath, amssymb, latexsym, amsthm}
\input epsf

\newcommand{\arx}[1]{\texttt{http://arxiv.org/abs/#1}}
\newcommand{\fA}{\mathfrak{A}}
\newcommand{\fB}{\mathfrak{B}}
\newcommand{\seq}[1]{\{#1\}_{n\in\N}}

\newcommand{\op}{\operatorname}
\newcommand{\maxfin}{\op{maxfin}}

\newcommand{\Cal}{\mathcal}

\newcommand{\B}{{\Cal B}}

\newcommand{\BG}{\B_\Gamma}

\newcommand{\BO}{\B_\Omega}

\newcommand{\F}{{\Cal F}}
\newcommand{\J}{{\Cal J}}

\newcommand{\N}{\naturals}
\newcommand{\NN}{{{}^{\naturals}\naturals}}

\renewcommand{\O}{\Cal O}

\newcommand{\R}{\reals}

\newcommand{\UU}{{\Cal U}}
\newcommand{\U}{\bigcup}
\newcommand{\V}{{\Cal V}}

\newcommand{\Impl}{\Rightarrow}

\long\def\forget#1\forgotten{}

\renewcommand{\b}{{\mathfrak b}}

\renewcommand{\d}{{\mathfrak d}}

\newcommand{\g}{\gamma}
\renewcommand{\i}{\item}

\newcommand{\w}{\omega}

\newcommand{\Iff}{\Leftrightarrow}

\newcommand{\nin}{\not\in}

\newcommand{\sbst}{\subseteq}

\newcommand{\sm}{\setminus}

\renewcommand{\pi}{pseudo-intersection}

\renewcommand{\>}{\rangle}
\newcommand{\E}{\exists}

\newcommand{\add}{{\sf add}}

\newcommand{\non}{{\sf non}}

\renewcommand{\t}{\tilde}

\newtheorem{thm}{Theorem}[section]
\newtheorem{prop}[thm]{Proposition}
\newtheorem{prob}[thm]{Problem}
\newtheorem{lem}[thm]{Lemma}
\newtheorem{cor}[thm]{Corollary}

\theoremstyle{definition}
\newtheorem{definition}[thm]{Definition}

\theoremstyle{remark}

\newcommand{\be}{\begin{enumerate}}
\newcommand{\ee}{\end{enumerate}}
\newcommand{\bi}{\begin{itemize}}
\newcommand{\ei}{\end{itemize}}

\newcommand{\sone}{{\sf S}_1}    \newcommand{\sfin}{{\sf S}_{fin}}
\newcommand{\ufin}{{\sf U}_{fin}}

\newcommand{\naturals}{{\mathbb N}}
\newcommand{\reals}{{\mathbb R}}

\author{Boaz Tsaban}
\address{Department of Mathematics and Computer Science, Bar-Ilan University,
Ramat-Gan 52900, Israel}
\email{tsaban@macs.biu.ac.il, http://www.cs.biu.ac.il/\~{}tsaban}

\title[Diagonalization between Hurewicz and Menger]
{A diagonalization property between Hurewicz and Menger}

\subjclass{%
Primary: 37F20; 
Secondary 26A03, 
03E75 
}

\keywords{%
Menger property, Hurewicz property, selection principles, continuous images%
}

\begin{document}
\begin{abstract}
In classical works, Hurewicz and Menger introduced
two diagonalization properties for sequences of open covers.
Hurewicz found a combinatorial characterization of these notions in
terms of continuous images.
Recently, Scheepers has shown that these notions are particular
cases in a large family of diagonalization schemas.
One of the members of this family is weaker
than the Hurewicz property and stronger than the Menger property,
and it was left open whether it can be characterized combinatorially
in terms of continuous images. We give a positive answer.
This paper can serve as an exposition of this fascinating subject.
\end{abstract}

\maketitle

\section{Introduction}
The following property was introduced by Menger \cite{MENGER}.
\begin{definition}
A set of reals $X$ has the \emph{Menger} property if for each
sequence $\seq{\UU_n}$ of
open covers
of $X$ there exist finite subsets $\F_n$ of $\UU_n$, $n\in\N$,
such that the collection $\seq{\cup\F_n}$ is a cover of $X$.
\end{definition}
For example, every compact, or even $\sigma$-compact, set of reals
satisfies the Menger property.

\begin{figure}[!h]
\begin{center}
\epsfysize=8 truecm {\epsfbox{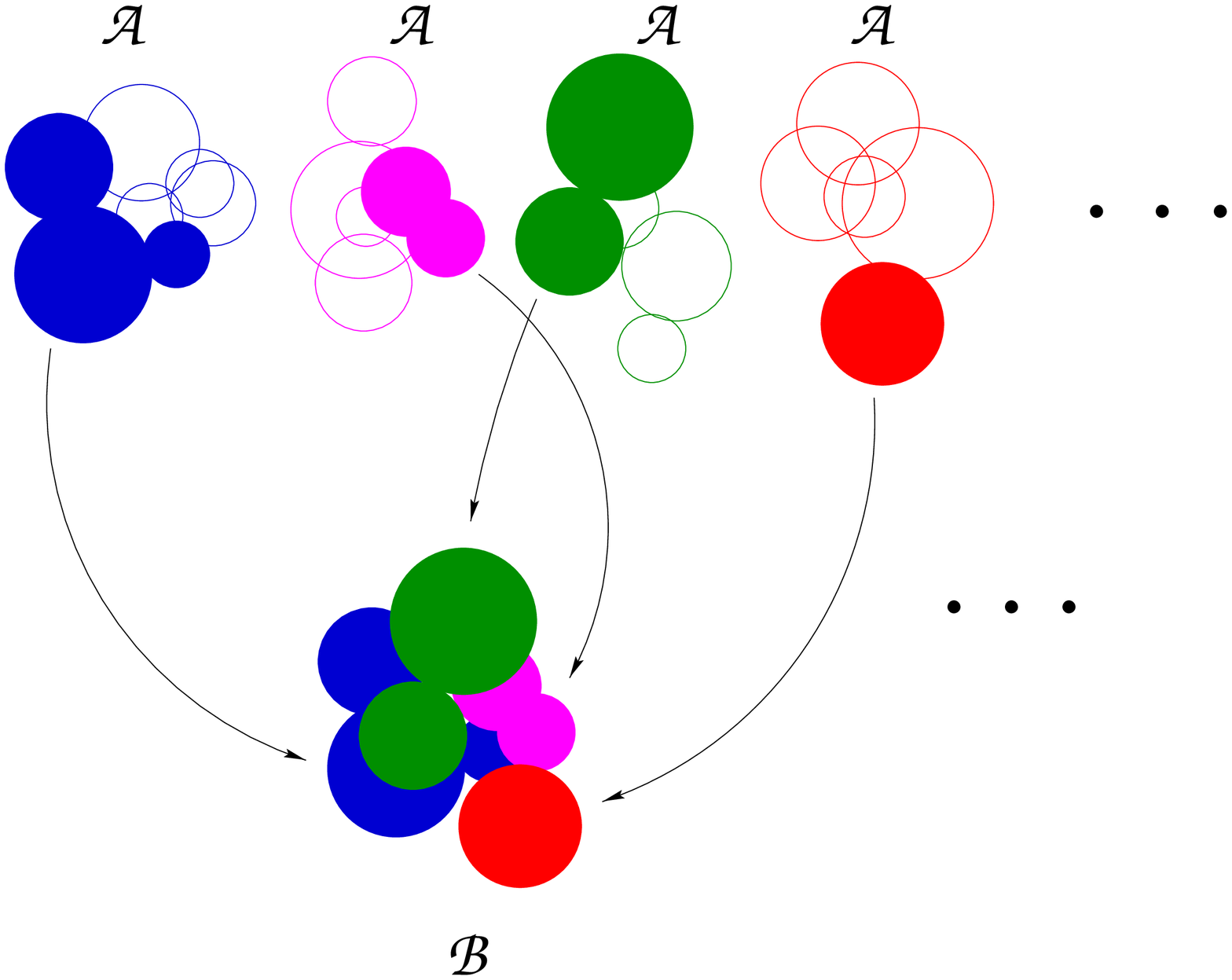}}
\end{center}
\caption{The selection principle $\ufin(\fA,\fB)$}\label{ufin}
\end{figure}

Following Hurewicz \cite{HURE25}, Scheepers introduced the following
diagonalization procedure of a sequence of covers \cite{comb1}
(see Figure \ref{ufin}):
\begin{definition}\label{ufindef}
Let $\fA$ and $\fB$ be collections of open covers of a space $X$.
$X$ has \emph{property $\ufin(\fA,\fB)$} if
for each sequence $\seq{\UU_n}$ of members of $\fA$
there exist finite subsets $\F_n$ of $\UU_n$, $n\in\N$,
such that either $\cup \F_n = X$ for some $n$, or else
the collection $\seq{\cup\F_n}\in\fB$.
\end{definition}

We will consider sets of reals for which the usual induced topology
has a subbase whose elements are \emph{clopen} (both closed and open),
that is, sets which are \emph{zero-dimensional}.
For convenience, we will also consider
the Baire space $\NN$ of infinite sequences
of natural numbers (equipped with the product topology).
The Baire space, as well as any separable and
zero-dimensional metric space
is homeomorphic to a set of reals; thus our results about sets of reals
are actually results about this more general case.

Let $X$ be a set of reals.
An $\w$-cover of $X$ is a cover such that each finite subset of
$X$ is contained in
some member of the cover. It is a $\g$-cover if it is infinite, and
each element of $X$ belongs to all but finitely many
members of the cover. Let $\O$, $\Omega$, and $\Gamma$
denote the collections of countable\footnote{%
There is no loss of generality in restricting attention to
countable covers here, since the spaces in question are
Lindel\"of.}
open covers, $\w$-covers, and
$\g$-covers of $X$, respectively.
In \cite{HURE27} Hurewicz studied the classes
$\ufin(\O,\O)$ (Menger property) and $\ufin(\O,\Gamma)$
(\emph{Hurewicz property}).
These diagonalization principles as well as several other natural
diagonalization principles (see Definition \ref{sonesfin}) were
studied by Scheepers, et.\ al., in a long series of papers
(\cite{comb1}, \cite{comb2}, \cite{comb3}, etc.).

For each of the diagonalization properties, it is desirable to have a
simple description of its underlying combinatorial structure.
For some of the properties this involves the quasiordering $\le^*$
defined on the Baire space $\NN$ by eventual
dominance:
$$f\le^* g\quad\mbox{ if}\quad f(n)\le g(n)\mbox{ for all but finitely many }n.$$
A subset $\F$ of $\NN$ is \emph{dominating} if for each $g$ in $\NN$ there
exists $f\in\F$ such that $g\le^* f$.

Hurewicz (\cite{HURE27}, see also Rec\l{}aw \cite{RECLAW})
has found the following elegant characterizations of the Menger and Hurewicz properties.
\begin{thm}[Hurewicz]
Let $X$ be a zero-dimensional set of reals.
\be
\i $X$ satisfies $\ufin(\O,\O)$ if, and only if,
every continuous image of $X$ in $\NN$ is not dominating.
\i $X$ satisfies $\ufin(\O,\Gamma)$ if,
and only if,
every continuous image of $X$ in $\NN$ is bounded (with respect to $\le^*$).
\ee
\end{thm}
Rec\l{}aw proved similar results for other important classes; see \cite{RECLAW}.

\section{The new property $\ufin(\O,\Omega)$}

Using finite unions, one can turn any countable cover into a $\g$-cover.
Thus, for each collection of covers $\fB$, the properties
$\ufin(\O,\fB)$, $\ufin(\Omega,\fB)$,
and $\ufin(\Gamma,\fB)$ are equivalent \cite{comb2}.
Therefore, the only (possibly) new property introduced by the general
scheme of Definition \ref{ufindef}
(when restricting attention to $\O$, $\Omega$, and $\Gamma$)%
\footnote{See \cite{tautau} for another type of covers which can
be smoothly incorporated into this framework.}
is $\ufin(\O,\Omega)$.
The property $\ufin(\O,\Omega)$ is weaker than the
Hurewicz property $\ufin(\O,\Gamma)$
and stronger than the
Menger property $\ufin(\O,\O)$, and according to \cite{comb2} it is
indeed new, that is, it is not provably equivalent to any of the
classical properties of Menger and Hurewicz.
Unlike these classical properties,
the combinatorial counterpart of the new property was less evident.
In \cite{CBC} it was proved that an analogous property 
(involving \emph{Borel} instead of open covers) can be characterized in
terms of the combinatorial structure of Borel images.
In Remark 10 of \cite{CBC} it was left open whether
a similar result can be obtained for $\ufin(\O,\Omega)$.
We answer this question positively. The proof is essentially an application of
Rec\l{}aw's arguments from \cite{RECLAW} to the proof of the
corresponding Borel result from \cite{CBC}.

For a finite subset $F$ of $\NN$, define $\max(F)\in\NN$ to be the
function $g$ such that $g(n)=\max\{f(n) : f\in F\}$ for each $n$.
In \cite{CBC} the following notion was introduced.
For a subset $Y$ of $\NN$,
$$\maxfin(Y):=\{\max(F) : F\mbox{ is a finite subset of }Y\}.$$
We prove the following characterization of the new property
$\ufin(\O,\Omega)$.
\begin{thm}\label{solved}
For a zero-dimensional set $X$ of reals, the following are equivalent:
\be
\i $X$ satisfies $\ufin(\O,\Omega)$;
\i For each continuous function $\Psi$ from $X$ to $\NN$, $\maxfin(\Psi[X])$ is not dominating.
\ee
\end{thm}
\begin{proof}
$2\Impl 1$: 
Assume that $\UU_n$, $n\in\N$, are open covers of $X$.
For each $n$, replacing each open member of $\UU_n$
with all of its clopen subsets we may assume that all elements of $\UU_n$
are clopen, and thus we may assume further that they are disjoint.
For each $n$ enumerate $\UU_n = \{U^n_m\}_{m\in\N}$.
As we assume that the elements $U^n_m$, $m\in\N$, are disjoint,
we can define a function $\Psi$ from $X$ to $\NN$ by
$$\Psi(x)(n)=m \Iff x\in U^n_m.$$
Then $\Psi$ is continuous. Therefore, $\maxfin(\Psi[X])$ is not dominating.
Let $g$ be a witness for that.
Then for each finite $F\sbst X$, $\max(\Psi[F])(n)<g(n)$ (i.e., $F\sbst\U_{k<g(n)}U^n_k$)
for infinitely many $n$.
Thus, $\seq{\U_{k<g(n)}U^n_k}$ is an $\w$-cover of $X$.

$1\Impl 2$: This was proved in \cite{CBC}. For completeness, we
give a minor variation of the original proof.
Since $\Psi$ is continuous, $Y=\Psi[X]$  also satisfies
$\ufin(\O,\Omega)$ \cite{comb2}.
Consider the basic open covers $\UU_n = \{U^n_m\}_{m\in\N}$ defined by
$U^n_m = \{f : f(n) = m\}$. Then there exist finite $\F_n\sbst\UU_n$, $n\in\N$,
such that either $Y=\cup\F_n$ for some $n$, or else
$\V=\{\cup\F_n : n\in\N\}$ is an $\w$-cover of $Y$.

The first case can be split into two sub-cases: If $Y=\cup\F_n$
for infinitely many $n$, then for these infinitely many $n$, the set $\{f(n) : f\in Y\}$ is
finite. Thus $\maxfin(Y)$ cannot be dominating.
Otherwise $Y=\cup\F_n$ for only finitely many $n$, therefore
we may replace each $\F_n$ satisfying $Y=\cup\F_n$ with $\F_n=\emptyset$,
so we are in the second case.
In the second case, since $Y\nin\V$ and $\V$ is an
$\w$-cover of $Y$, we have that each
finite subset of $Y$ is contained in infinitely many elements of $\V$.
Define $g\in\NN$ by $g(n)=\max\{m : U^n_m\in\F_n\}+1$.
For each finite $F\sbst Y$, we have that
$F\sbst\cup\F_n$ and thus $\max(F)(n) < g(n)$ for infinitely many $n$. Then $g$
witnesses that $\maxfin(Y)$ is not dominating.
\end{proof}

The following diagonalization properties, which generalize some other
classical notions, were introduced by Scheepers \cite{comb1}.
\begin{definition}\label{sonesfin}
For a set $X$ of reals, define the following diagonalization properties:
\bi
\item[$\sone(\fA,\B)$:]
{For each sequence $\seq{\UU_n}$ of members of $\fA$, there is a sequence
$\seq{U_n}$ such that for each $n$ $U_n\in \UU_n$, and
$\seq{U_n}\in\fB$.}
\item[$\sfin(\fA,\fB)$:]
{For each sequence $\seq{\UU_n}$ of members of $\fA$, there is a sequence
$\seq{\F_n}$ such that each $\F_n$ is a finite subset of $\UU_n$, and
$\U_{n\in\N}\F_n\in\fB$.}
\ei
\end{definition}

Let $\B$, $\BO$, and $\BG$
denote the collections of countable Borel covers, $\w$-covers, and
$\g$-covers of $X$, respectively.
Using Theorem \ref{solved} and results from \cite{CBC},
we get the following corollary, which answers another question from \cite{CBC}.
(Recall that $\ufin(\O,\Omega)=\ufin(\Gamma,\Omega)$.)
\begin{cor}
For a zero-dimensional set $X$ of reals, the following are equivalent:
\be
\i $X$ satisfies $\sone(\BG,\BO)$,
\i Every Borel image of $X$ satisfies $\sone(\Gamma,\Omega)$,
\i Every Borel image of $X$ satisfies $\sfin(\Gamma,\Omega)$;
\i Every Borel image of $X$ satisfies $\ufin(\Gamma,\Omega)$.
\ee
\end{cor}
\begin{proof}
This follows from Theorem \ref{solved}, together with the equivalences
$\sone(\BG,\BO)=\sfin(\BG,\BO)=\ufin(\BG,\BO)$, and the fact that
$\sone(\BG,\BO)$ is closed under taking Borel images, see \cite{CBC}.
\end{proof}

\section{Finite dominance}

For a subset $Y$ of $\NN$, the property that $\maxfin(Y)$ is not
dominating can be stated in terms of finite dominance, a notion
defined by Blass in a recent study \cite{Andreas}.

\begin{definition}
Let $Y$ be a subset of $\NN$, and $k\in\N$. $Y$ is \emph{$k$-dominating}
if for each $f\in\NN$ there exist $g_1,\ldots, g_k\in Y$
such that $f(n)\le \max\{g_1(n),\ldots,g_k(n)\}$ for all but
finitely many $n$.
\end{definition}
In other words, $Y$ is $k$-dominating when the collection
$\{\max(F) : F\sbst Y\mbox{ and }|F|=k\}$ is dominating.

\begin{prop}\label{andreas}
For $Y\sbst\NN$, $\maxfin(Y)$ is dominating if, and only if,
there exists a natural number $k$ such that $Y$ is $k$-dominating.
\end{prop}
\begin{proof}
We will prove the less trivial implication.
Assume that for each $k$, $Y$ is not $k$-dominating.
For each $k$, let $g_k\in\NN$ witness that $Y$ is not
$k$-dominating. Since the collection $\{g_k : k\in\N\}$
is countable, there exists $f\in\NN$ bounding it with
respect to eventual dominance $\le^*$.
Then $f$ witnesses that $\maxfin(Y)$ is not dominating.
\end{proof}

\begin{cor}
For a zero-dimensional set $X$ of reals, the following are equivalent:
\be
\i $X$ satisfies $\ufin(\O,\Omega)$;
\i For each continuous function $\Psi$ from $X$ to $\NN$
and for each natural number $k$, $\Psi[X]$ is not $k$-dominating.
\ee
\end{cor}

\section{Reduced products}

The characterization in Theorem \ref{solved} has an elegant
statement in the language of filters.
Let $\F$ be a filter over $\naturals$. An equivalence
relation $\sim_\F$ is defined on $^\naturals\naturals$ by
$$f\sim_\F g \Leftrightarrow \{n : f(n) = g(n)\} \in \F.$$
The equivalence class of $f$ is denoted $[f]_\F$, and the set
of these equivalence classes is denoted $^\naturals\naturals/\F$.
Using this terminology, $[f]_\F < [g]_\F$ means
$$\{n: f(n) < g(n)\} \in \F.$$

\begin{lem}[\cite{CBC}]\label{maxfinlemma}
Let $Y\sbst\NN$ be such that for each $n$ the set
$\{h(n):h\in Y\}$ is infinite. Then the following are equivalent:
\begin{enumerate}
\item{$\maxfin(Y)$ is not a dominating family.}
\item{There is a non-principal filter $\F$ on $\naturals$ such that the subset $\{[f]_{\F}:f\in Y\}$ of the reduced product $^\naturals\naturals/\F$ is bounded.}
\end{enumerate}
\end{lem}

\begin{cor}\label{maxfin}
For a zero-dimensional set $X$ of reals, the following are equivalent:
\be
\i $X$ satisfies $\ufin(\O,\Omega)$;
\i For each continuous function $\Psi$ from $X$ to $^{\naturals}\naturals$,
either there is a principal filter
${\Cal G}$ for which $\Psi[X]/{\Cal G}$ is finite, or else there is a
nonprincipal filter ${\Cal F}$ on $\naturals$
such that the subset $\Psi[X]/{\Cal F}$ of the reduced product
$^{\naturals}\naturals/{\Cal F}$ is bounded.
\ee
\end{cor}
\begin{proof}
This follows from Theorem \ref{solved}
and Lemma \ref{maxfinlemma}, see \cite{CBC} for
the proof of the corresponding Borel version.
\end{proof}

\forget
\section{Additivity numbers}

The \emph{critical cardinality} of a class $\J\subset P(\R)$, $\non(\J)$,
is the minimal cardinality of a set of reals $X$ such that $X\nin\J$.

Let $\b$ be the minimal cardinality of an unbounded family
in $\NN$, and $\d$ be the minimal cardinality of
a dominating family in $\NN$ \cite{vD}.
From Theorem \ref{solved} one immediately gets the following.
\begin{cor}
$\non(\ufin(\O,\Omega))=\d$.
\end{cor}
\begin{proof}
If $|Y|<\d$, then $|\maxfin(Y)|<\d$ as well.
On the other hand, from the definition of $\d$, there exists
a dominating family of cardinality $\d$.
\end{proof}
This result was first proved in \cite{comb2}.

The \emph{additivity number} of a class $\J$ with $\cup\J\nin\J$,
$\add(\J)$, is the minimal cardinality of a collection $\F\sbst\J$ such that
$\cup\F\nin\J$.
\begin{lem}\label{combadd}
\be
\i Assume that $\kappa<\b$, and that $X_\alpha$, $\alpha<\kappa$,
are bounded subsets of $\NN$.
Then the set $X=\U_{\alpha<\kappa}X_\alpha$ is bounded.
\i Assume that $\kappa<\b$, and that $X_\alpha$, $\alpha<\kappa$,
are subsets of $\NN$ which are not dominating.
Then $X=\U_{\alpha<\kappa}X_\alpha$ is not dominating.
\ee
\end{lem}
\begin{proof}
(1) For each $\alpha<\kappa$, let
$g_\alpha$ bound $X_\alpha$.
Let $G=\{g_\alpha : \alpha<\kappa\}$.
As $\kappa<\b$, $G$
is bounded by some function $h$.
Then $h$ bounds $X$.

(2) For each $\alpha<\kappa$, let
$g_\alpha$ witness that $X_\alpha$ is not dominating.
Let $G=\{g_\alpha : \alpha<\kappa\}$.
As $\kappa<\b$, $G$ is bounded by some function $h$.
It is easy to see that $h$ witnesses that $X$ is not a dominating family.
\end{proof}

In \cite{comb2} it is pointed out that the
Menger and Hurewicz properties
are closed under taking countable unions.
This fact can be derived from the definition of these classes, by
partitioning the given sequence of covers to infinitely many infinite
sequences. Using the characterizations of these properties in terms
of their images in $\NN$, we can get more information on
their additivity numbers.

\begin{thm}\label{adds}
\be
\i $\add(\ufin(\B,\BG))=\add(\ufin(\O,\Gamma))=\b$,
\i $\b\le\add(\ufin(\B,\B))=\add(\ufin(\O,\O))\le\d$;
\ee
\end{thm}
\begin{proof}
This follows from Lemma \ref{combadd}
and the characterizations of these classes
in terms of Borel \cite{CBC} and continuous images, together
with the fact that
$\add(\J)\le\non(\J)$ whenever all singletons are members of $\J$.
\end{proof}

\begin{prob}
Is $\add(\ufin(\O,\O))$ provably equal to $\b$ or $\d$?
\end{prob}

On Problem 5 of \cite{comb2} it is asked, in particular, whether
the property $\ufin(\O,\Omega)$ is closed under taking finite or
countable unions. The simplest way to obtain a positive result
would be to prove the corresponding combinatorial version
as in Lemma \ref{combadd}.
Unfortunately this is wrong.
\begin{prop}\label{not-working1}
There exist subsets $A$ and $B$ of $\NN$ such that
$\maxfin(A)$ and $\maxfin(B)$ are not dominating,
but $\maxfin(A\cup B)$ is dominating.
\end{prop}
\begin{proof}
Take any dominating set $D$, and consider the following sets:
Let $E$ be the set of even natural numbers, and $O$ be the
set of odd natural numbers. Choose
\begin{eqnarray}\nonumber
A & = & \{f : (\E g\in D)(\forall \mbox{ even } n)\ f(n)=g(n)
\mbox{ and } f\restriction O\equiv 0 
\}\cr\nonumber
B & = & \{f : (\E g\in D)(\forall \mbox{ odd } n)\ f(n)=g(n)
\mbox{ and } 
f\restriction E\equiv 0\}\cr\nonumber
\end{eqnarray}
\end{proof}

Proposition \ref{not-working1} does not imply that $\ufin(\O,\Omega)$ is
not closed under taking finite unions:
For a subset $X$ of $\NN$ and a strictly increasing $h\in\NN$, define
$X\circ h = \{f\circ h : f\in A\}$. Then $X\circ h$ is a continuous
image of $X$. Now, for $f(n)=2n$ and $g(n)=2n+1$, we have that
$A\circ f$ and $B\circ g$ are dominating, thus $A$ and $B$ do not
satisfy $\ufin(\O,\Omega)$.

The next step would be to consider unions of sets $X$ such that
for any $h\in\NN$, $X\circ h$ is not dominating. In answer to my question,
Gitik has proved that this stronger assumption is also not enough.

\begin{thm}[Gitik]\label{not-working2}
Assume the Continuum Hypothesis.
There exist subsets $A$ and $B$ of $\NN$ such that
for any $h\in\NN$, $\maxfin(A\circ h)$ and $\maxfin(B\circ h)$
are not dominating, but $\maxfin(A\cup B)$ is dominating.
\end{thm}
\begin{proof}
The proof is a little technical and the reader unfamiliar with
this sort of arguments may wish to skip it.

Let $\<X_\alpha : \alpha<\aleph_1\>$ be an enumeration of all infinite subsets of $\N$.
Let $S_1,S_2,\dots$ be a partition of $\N$ into infinitely many infinite sets.
Define $\<A_\alpha : \alpha<\aleph_1\>$ by induction.
Assume that that $A_\beta$ is defined for all $\beta<\alpha$.
Define $A_\alpha$ as follows.
Let $\<Y_\alpha : \alpha\in\N\>$ be an enumeration of
all infinite sets which are finite Boolean combinations of sets from
$\{X_\beta : \beta\le \alpha\}\cup \{A_\beta : \beta<\alpha\}$.
Define, by induction on $k$, elements $a_k\in A_\alpha$ and $b_k\in \N\sm A_\alpha$:
For each $k$, Let $n_k$ be such that $k\in S_{n_k}$.
Then we let $a_k$ and $b_k$ be any
two distinct elements of $Y_{n_k} \sm \{a_l, b_l : l<k\}$.

Let $D=\<d_\alpha : \alpha<\aleph_1\>$ be a dominating family in $\NN$
such that for $\alpha<\beta$, $d_\alpha\le^* d_\beta$.
For each $\alpha<\aleph_1$ let $f_\alpha$ be equal to $d_\alpha$ on $A_\alpha$ and to $0$ on $\N\sm A_\alpha$,
and $g_\alpha$ be $d_\alpha$ on $\N\sm A_\alpha$ and $0$ on $A_\alpha$.
Set $A = \{f_\alpha : \alpha<\aleph_1\}$ and $B = \{g_\alpha : \alpha<\aleph_1\}$.

Towards a contradiction, assume that for some $h\in\NN$, $\maxfin(A\circ h)$ is dominating.
Find $\beta<\aleph_1$ such that $h[\N]=X_\beta$.
As the family $\maxfin\{f_\alpha\circ h : i\le \beta\}$ is countable,
it is not dominating; let $u\in\NN$ be a witness for that.
We may assume that $n<u(n)$ for all $n$.
By the assumption there is a finite subset $F$ of $\aleph_1$ such
that
$s=\max\{f_\alpha\circ h: \alpha\in F\}$ dominates $u$.
By the choice of $u$ we can assume that $F$
is a subset
of $\aleph_1\sm\beta$. But then $X_\beta\sm\cup\{A_\alpha : \alpha\in F\}$ is infinite and the value
of $s$ on this set is $0$.
Thus $s$ cannot dominate $u$.

A similar argument shows that there exists no $h\in\NN$ such that
$\maxfin(B\circ h)$ is dominating. But clearly $\maxfin(A\cup B)$ is
dominating.
\end{proof}

Thus, the best we can obtain using this approach is the following result.
For classes $\J$ and $\t\J$ with $\cup\J\nin\t\J$,
we write $\add(\J,\t\J)$ for the minimal cardinality of a collection
$\F\sbst\J$ with $\cup\F\nin\t\J$ \cite{jubar}.

\begin{lem}
Assume that $\kappa<\d$, and that $X_\alpha$, $\alpha<\kappa$,
are bounded subsets of $\NN$.
Then for $X=\U_{\alpha<\kappa}X_\alpha$, $\maxfin(X)$ is not dominating.
\end{lem}
\begin{proof}
For each $\alpha<\kappa$, let
$g_\alpha$ bound $\maxfin(\Psi[X_\alpha])$.
Let $G=\{g_\alpha : \alpha<\kappa\}$.
As $\kappa<\d$, the family $\maxfin(G)$
is not dominating -- let $h$ witness this.
Assume that $F$ is a finite subset of $X$.
Then for each $f\in F$ there exists $g_f\in G$ such that
$f\le^* g_f$.
For infinitely many $n$, $g_f(n)\le h(n)$ for all $f\in F$.
Thus, for infinitely many $n$, $f(n)\le h(n)$ for all $f\in F$.
\end{proof}

\begin{cor}
$\add(\ufin(\B,\BG),\ufin(\B,\BO))=\add(\ufin(\O,\Gamma),\ufin(\O,\Omega))=\d$.
\end{cor}
Thus, a union of less than $\d$ many sets having the Menger property
must satisfy the new property $\ufin(\O,\Omega)$.
\forgotten

\textbf{Acknowledgments.}
I wish to thank Ireneusz Rec\l{}aw for reading this paper and
detecting some typos. I also thank Andreas Blass for the observation
\ref{andreas}.

\end{document}